\documentclass{article}
\usepackage{amssymb,amsmath,amsthm,xypic,graphicx,marvosym,xy}
\usepackage{color,epsfig}
\usepackage{newlfont}
\theoremstyle{plain}

\newcommand\Mod{\mathfrak{M}}

\newcommand\C{\mathbb{C}}

\newcommand\M{\overline{\mathfrak{M}}}

\newcommand\Mn{\Mod_{0,n}}

\newtheorem*{thm*}{Theorem}
\newtheorem{thm}{Theorem}[section]
\newtheorem{defn}[thm]{Definition}

\newtheorem{ex}[thm]{Example}

\newtheorem{prop}[thm]{Proposition}

\newtheorem{notate}[thm]{Notation}
\title{A polygonal presentation of $Pic(\M_{0,n})$ }
\author{Sarah Carr}
\date{}
\begin{document} 
\maketitle
\begin{abstract}
This article gives a combinatorial presentation
of the Picard group of $\M_{0,n}$.  In the first
section, we recall Keel's well-known presentation of $Pic(\M_{0,n})$
using boundary divisors of $\M_{0,n}$ as generators,
and describe a basis for $Pic(\M_{0,n})$ recently discovered by Keel
and Gibney.  In 
the second
section, we present a theorem
which gives an explicit and very simple expression for every boundary divisor
in terms of this basis, thereby yielding a new presentation with a
minimal set of relations.
\end{abstract}

\section{The Picard group of $\M_{0,n}$}

We begin by presenting definitions and theorems
about divisors on
$\M_{0,n}$, the
compactification of the moduli space of genus 0 $n$-pointed curves over $\C$.

\begin{defn} Let $X$ be a smooth manifold, and let $Div(X)$ be the
  group formally generated by Weil divisors on $X$.  The Picard group,
  $Pic(X)$, 
  is the quotient of $Div(X)$ by the principal divisors.\end{defn}

We have the following
characterization/definition of the Picard group of Weil divisors on
$\M_{0,n}$.

\begin{thm}\cite{Ke} The Picard group, $Pic(\M_{0,n})$, is isomorphic
  to the quotient
  $Div(\M_{0,n})/\sim$, where $\sim$ denotes numerical equivalence of
  divisors.
\end{thm}

$\M_{0,n}$ denotes the
Deligne-Mumford compactification of
$\Mn$ \cite{DM}.  The points of $\Mn$ are Riemann
spheres with $n$ marked points
modulo isomorphism.  Let $S_n=\{1,...,n\}$ denote the alphabet used to index the set of marked points, $Z_n=\{z_1,...,z_n\}$, on $\Mn$.
 The boundary divisor, $\M_{0,n}\setminus \Mn$, parametrizes
stable curves of type ($0,n$). The stable curves are curves which are nodal Riemann surfaces
of genus zero such that each component has at least three marked or singular
points.  In this article we will use the term boundary divisor to
refer to the irreducible components of $\M_{0,n} \setminus \Mn$.

\begin{figure}
\begin{center}
 \scalebox{.8}{\input{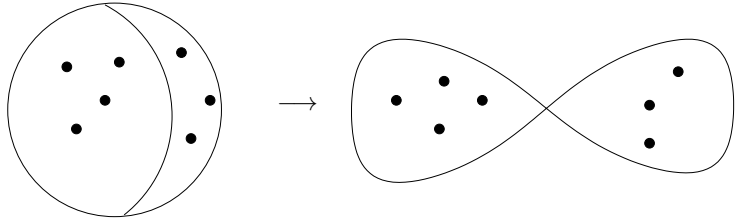}}
\caption{A point on a boundary divisor in $\M_{0,n}$}
\label{pinch}
\end{center}
\end{figure}

The (irreducible) boundary divisors can be combinatorially enumerated by specifying
a partition of $S_n$ into two subsets, $A$ and $S_n\setminus A$, with
$2\leq |A| \leq n-2$ in the following way.  Any simple closed loop on a sphere
with $n$ marked points partitions the points of $Z_n$ (indexed by $S_n$) into two subsets
as in the left hand side of figure \ref{pinch}.  
Pinching this loop to a point yields a nodal topological surface as in the right hand side of
figure \ref{pinch}.  The stable curves of this type are obtained by
putting all possible complex structures on this topological surface.
A single boundary component parametrizes these stable curves for a
given pinched loop.

We
denote by $d_A$
the boundary divisor in which the loop pinches the subset indexed by
$A\subset S_n$, hence $d_A=d_{S_n\setminus A}$.  We denote the set of
(irreducible) boundary divisors on $\M_{0,n}$ by $D^n$.  This set has
cardinality $2^{n-1}-1-n$.

\begin{ex}

The set of boundary divisors, $D^4$, on $\M_{0,4}$ contains the
three divisors, $d_{1,2}, d_{1,3}$ and $d_{1,4}$.

The set of boundary divisors, $D^5$, on $\M_{0,5}$ contains the 10
divisors, $d_A$, where $A\subset S_5$ has cardinality 2. 
\end{ex}

We now introduce the geometric notation for a boundary divisor which will be useful in the proofs in this paper.   We will frequently refer to cyclic orderings, denoted $(i_1,...i_n)$, on the set $S_n$.  One may think of a cyclic ordering as an $n$-sided, oriented polygon, whose vertices are labeled by the elements of $S_n$ such that $i_{p+1}$ follows $i_p$ in the clockwise direction, and hence $(i_1,...,i_n)=(i_2,...,i_n,i_1)= \cdots =(i_n, i_1..., i_{n-1})$. 
We denote by $P_n$ the $n$-sided, oriented polygon whose cyclic ordering is $(1,...,n)$ and call this the standard ordering.

\begin{notate}\label{divnot} To obtain a geometric representation of the
boundary divisor $d_A$ associated to $A=\{i_1,...,i_p\}$, we partition $A$ into
subsets, $B_1,...B_k$, of indices which are adjacent on the polygon $P_n$, with
the condition that there is a unique minimal $k$.  In other words, there is no
union of two sets, $B_a$ and $B_b$ that form an adjacent block on $P_n$.  Then,
the set $\{B_1,...,B_k\}$ inherits a cyclic ordering from $P_n$.  Without loss
of generality, we may assume that this cyclic ordering corresponds to the
subscripts on the $B_i$ and that the number 1 is an element both of $B_1$ and
of $A$.  Given these conditions, we have that for each $B_i$ there exists a gap,
$G_i$, of adjacent indices on $P_n$ which are not in $A$ and which are exactly
the indices between $B_i$ and $B_{i+1}$ on $P_n$. Then we denote the divisor
$d_A$ by the $2k$-sided polygon whose vertices are cyclically labeled by $(B_1,
G_1, B_2, G_2, ... ,G_{k-1}, B_k, G_k)$.  In the case where $k=1$, we draw the
divisor polygon as circle with two points, $B_1$ and $G_1$.
 \end{notate}

\begin{ex}
For $n=10$ two divisors are represented in figure 2.  The divisor $d_{\{1,2,4,7,10\}}$ is represented by a hexagon and the divisor $d_{\{1,2,3,9,10\}}$ is represented by a two-pointed circle.

\begin{figure}
\begin{center}
\scalebox{.6}
{\input{pillowgon5.pstex_t}}
\caption{Boundary divisors in $\M_{0,10}$}
\label{pillowgon5}
\end{center}
\end{figure}
\end{ex}

\begin{thm}\cite{Ke}\label{kepres}
A presentation of 
the Picard group, $Pic(\M_{0,n})$, is given by taking 
the classes, $\delta_A$ of the boundary divisors, $d_A\in D^n$ as generators, subject to the
following relations: for any four distinct elements,
$i,j,k,l$ in $S_n$,
$$\sum_{\textstyle{{i,j \in A}\atop{k,l\notin A}}} \delta_A =
\sum_{\textstyle{{i,k \in A}\atop{j,l\notin A}}} \delta_A =
\sum_{\textstyle{{i,l \in A}\atop{j,k\notin A}}} \delta_A.$$
\end{thm}


The following proposition due to A. Gibney \cite{GM} specifies a basis
 for $Pic(\M_{0,n})$ and also yields an expression for its dimension.

\begin{prop}\label{Keel}
Let $({i_1},...,{i_n})$ denote any cyclic ordering.  Then a basis for $Pic(\M_{0,n})$ is given by the
divisors defined by
nonempty subsets of marked points on the $n$-gon which do not form an
adjacent set of vertices on the $n$-gon.  We call this set of divisors
the {\bf non-adjacent basis} for the cyclic ordering $(i_1,...,i_n)$. 
\end{prop}

In the notation \ref{divnot}, this proposition means that the set of
divisors $$\{(B_1,G_1,...,B_k,G_k):k\geq 2\}$$ is a basis for $Pic(\M_{0,n})$. 
Therefore the main purpose of this paper is to find an expression for the
divisors denoted $(B_1,G_1)$.

{\bf Remarks}.
The following combinatorial formula for the dimension follows
immediately from counting the elements of a non-adjacent basis:
\begin{equation}
\mathrm{dim}(Pic(\M_{0,n}))=2^{n-1}-1-n- {n\choose
  2}+n=2^{n-1}-1-{n\choose 2}.
\end{equation}
This dimension was calculated by S. Keel in \cite{Ke} as the dimension of the first Chow group of $\M_{0,n}$.
%

\begin{ex}
Consider the standard ordering,
$(1,2,3,4,5)$.  Then the non-adjacent basis for
$Pic(\M_{0,5})$ for this ordering is
given by the five divisor classes $$\delta_{\{1,3\}}, \delta_{\{ 1,4\}}, \delta_{\{2,4\}},
\delta_{\{2,5\}}, \delta_{\{3,5\}}.$$
\end{ex}
\begin{proof} (of Proposition \ref{Keel}). Without loss of generality, we prove this proposition for the standard cyclic ordering, $(1,2,...,n)$.  We restate this theorem by saying that the set of divisor classes, $\{ \delta_A$ such that by notation \ref{divnot}, $A=B_1\sqcup \cdots \sqcup B_k$ where $k\geq 2\}$, forms a basis for $Pic(\M_{0,n})$.  We will prove this proposition by showing that all divisors corresponding to $k=1$ may be written as divisors corresponding to $k\geq 2$ and then by justifying that this set of ``non-adjacent'' divisor classes indeed forms a basis.

Consider any divisor class, $\delta_I$ where $I$ is a set of adjacent indices taken modulo $n$, $\{m,m+1,...,m+p\}$.  By taking $[i,j,k,l]=[m-1,m,m+p,m+p+1]$, we have the following relation from theorem \ref{kepres}:
\begin{equation}\label{nabaseq}
\sum_{\textstyle{{m-1,m+p \in A}\atop{m,m+p+1\notin A}}} \delta_A =
\sum_{\textstyle{{m-1,m+p+1 \in A}\atop{m,m+p\notin A}}} \delta_A.
\end{equation}
The left-hand side of equation \eqref{nabaseq} is the sum of divisor classes, $\delta_A$ where $A$ is a set of non-adjacent indices, since the partition of any such $A$ is at least two sets of consecutive indices, one containing $m-1$ and the other containing $m+p$.  The sets containing $m-1$ and $m+p$ cannot be the same, because any consecutive ordering containing both $m-1$ and $m+p$ must contain either $m$ or $m+p+1$, but there are no such terms on the left-hand side.

 The right-hand hand side is the sum of divisor classes, $\delta_A$, such that $A$ contains $m-1$ and $m+p+1$,
 and such that both $m$ and $m+p$ are contained in gaps.  The only divisor class term indexed by a consecutive
 set 
is $\delta_A$, $A=S_n\setminus I$.  Therefore for any set of adjacent indices $I$, we may write equation \eqref{nabaseq} as
$$\sum_p \delta_{A_p} =\delta_I + \sum_q \delta_{B_q}$$
where the $A_p, B_q$ are all elements of the non-adjacent basis, and hence these divisor classes span the Picard group.

Furthermore, there are exactly $2^{n-1}-1-{n\choose 2}$ non-adjacent divisors, which is the dimension of the Picard group \cite{Ke}, so these elements indeed form a basis.
\end{proof}

\section{A polygonal presentation of $Pic(\M_{0,n})$}

In this section, we give a simple expression of each boundary divisor
in $Pic(\M_{0,n})$ in
terms of any non-adjacent basis.
This yields a new and very simple presentation for $Pic(\M_{0,n})$
with a minimal set of relations.


\begin{thm} Let $\gamma$ denote a cyclic ordering $({i_1},\ldots,
i_n)$ on $S_n$.  Then
$Pic(\M_{0,n})$ is generated by the set of classes, $\delta_I$, of boundary
divisors of
$\M_{0,n}$, subject to the relations

\begin{equation}\label{sarahsthm}\delta_I=\sum_{J\in {\cal J}}
  \delta_J - \sum_{K\in {\cal K}} 
\delta_K,
\end{equation}
where $I$ denotes a consecutive subset of points for the ordering
$\gamma$, ${\cal J}$ denotes the set of non-adjacent subsets 
$$J=B_1\cup \cdots \cup B_j$$ 
of $\{1,\ldots,n\}$ such that $I$ is equal to a ``segment'' of even
length
$B_i,G_i,\ldots,B_k,G_k$ or $G_i,B_{i+1},G_{i+1},\ldots,G_{k},B_{k+1}$
of $(B_1,G_1,\ldots,B_N,G_N)$, and ${\cal K}$ denotes the set of
non-adjacent subsets $K=B_1\cup \cdots\cup B_j$ such that
$I$ is equal to a ``segment'' of odd length 
$B_i,G_i,\ldots,B_k,G_k,B_{k+1}$ or
$G_i,B_{i+1},G_{i+1},\ldots,B_k,G_k$
of $(B_1,G_1,\ldots,B_N,G_N)$.
\end{thm} 

The beauty of the theorem is more easily seen by rephrasing it as: the coefficients of any divisor in the basis of the Picard group given by a cyclic ordering can be calculated by the parity of the defining blocks of the divisor.  
The precise statement of the theorem does not do justice to its simplicity, as illustrated in the following example.

\begin{ex}
We have the following expression for the divisor, $\delta_{\{1, 2, 3\}}$, in the basis of $Pic(\M_{0,6})$ given by the cyclic ordering $(1,2,3,4,5,6)$:
$$\delta_{\{1,2,3\}} = -\delta_{\{1,3\}} + \delta_{\{1,4\}} +\delta_{\{3,6\}}-\delta_{\{4,6\}}+\delta_{\{1,2,4\}}-\delta_{\{1,3,5\}}+\delta_{\{1,4,5\}}.$$
\end{ex}

\begin{proof}
Without loss of generality, we may prove this theorem on the standard cyclic ordering $\gamma =(1,2,...,n)$.
Let $\Delta_\gamma$ be set of divisors which is a basis for
$Pic(\M_{0,n})$ with respect to the cyclic order $\gamma$ by proposition
\ref{Keel}.  We denote by 
$\delta_{B_1\cdots B_N} = [B_1,G_1,\ldots,B_N,G_N]$ an element of
$\Delta_\gamma$.  Let $I=(m,...,m+p)$ (where indices are taken modulo 
$n$) be a consecutive subset for the cyclic order
$\gamma$.

Then we
may restate the theorem as follows.  One can express $\delta_I$ as the
linear combination of elements of $\Delta_\gamma$:
\begin{equation}\label{deltaI}
\delta_I = \sum C_{B_1,...,B_{J_k}}^I [B_1,G_1...,B_{J_k}, G_{J_k}],
\end{equation} and the coefficients are given by
\begin{equation}\label{piccoeffs}
C^I_{B_1,...,B_N} = \begin{cases} 
1 & I=\bigcup_{p=1}^j B_{i+p} \cup
  G_{i+p}\  \\
& \qquad\Bigl(\hbox{or } I=G_{i} \cup \bigl( \bigcup_{p=1}^j B_{i+p} \cup G_{i+p} \bigr) \cup
 B_{i+j+1} \Bigr)\\
-1 & I=\bigl(\bigcup_{p=1}^j B_{i+p} \cup
  G_{i+p} \bigr) \cup B_{i+j+1}\\\ & \qquad \Bigl(\hbox{or } I=G_i\cup \bigl(\bigcup_{p=1}^j B_{i+p} \cup
  G_{i+p} \bigr)\Bigr) \\
0 &\hbox{otherwise},
\end{cases}
\end{equation}
where for $1<p<j$, $i+p$ is taken modulo $n$.

To prove this theorem, we take $[i,j,k,l] = [m-1,m,m+p,m+p+1]$ as in equation \eqref{nabaseq},
\begin{equation}\label{nabaseq2}
\sum_{\textstyle{{m-1,m+p \in A}\atop{m,m+p+1\notin A}}} \delta_A =
\sum_{\textstyle{{m-1,m+p+1 \in A}\atop{m,m+p\notin A}}} \delta_A.
\end{equation}

I claim that the left-hand side of equation \eqref{nabaseq2} is a sum of all of the divisor classes, $\delta_A$, such that $A$ the $I$ is the disjoint union of an even number ($\geq 2$), of $B_p,G_q$ (as in notation \ref{divnot}) and the right-hand side of this equation is the sum of $\delta_A$ such that $I$ can be expressed as the disjoint union of an odd number ($\geq 1$) of $B_p, G_q$.  The justification of this claim proves the theorem.

Any divisor class, $\delta_A$, which is a term on the left-hand side of equation \eqref{nabaseq2} has a polygonal representation as in figure 3.  Therefore $I$ is the disjoint union of the even number of sets, $G_1\sqcup B_2 \sqcup G_2\sqcup \cdots G_{h-1} \sqcup B_h$.  While any divisor class, $\delta_A$, which is a term on the right-hand side of \eqref{nabaseq2} has a polygonal representation as in figure 3, so that $I$ is the disjoint union of an odd number of sets, $G_1\sqcup B_2\sqcup \cdots \sqcup B_h \sqcup G_h$.

\begin{figure}
\begin{center}
 \scalebox{.6}{\input{pillowgon6.pstex_t}}
\label{pillowgon6}
\end{center}
\end{figure}

By bringing all terms except $\delta_I$ over to the left-hand side of equation \eqref{nabaseq2}, we have the desired expression in the statement of the theorem.

\end{proof}
\vspace{1cm}
{\it Acknowledgements}.  This paper was written at the University of
Pennsylvania and is a product of discussions with A. Gibney,
D. Krashen and L. Schneps during my stay there.  I would like to thank
the University and these three researchers for their support and
clarifications.\\

\pagebreak 
\noindent Sarah Carr\\
Universit\'e Paris XI\\
Dept. de math\'ematiques B\^atiment 425\\
91405 Orsay Cedex FRANCE\\
sarah.carr@math.u-psud.fr

\end{document}